\documentclass[a4paper]{article}

\usepackage{amssymb, amsmath, amsthm, latexsym, verbatim}

\newcommand \K{\delta}
\newcommand \la{\lambda}
\newcommand \Rn{\mathbb R^n}
\newcommand \rnu{\mathbb R^{\nu}}
\newcommand \rk{\mathrm{rk}}
\newcommand \Ker{\mathrm{Ker}}
\renewcommand \Im{\mathrm{Im}}
\newcommand \Cliff{\mathrm{Cliff}}
\newcommand \Span{\mathrm{Span}}

\newcommand \Hom{\mathrm{Hom}}
\newcommand \SO{\mathrm{SO}}
\newcommand \diag{\mathrm{diag}}
\renewcommand \L{\Lambda}
\newcommand \T{\theta}
\renewcommand \a{\alpha}
\renewcommand \b{\beta}
\renewcommand \O{\mathrm{O}}
\newcommand \<{\langle}
\renewcommand \>{\rangle}

\newtheorem{theorem}{Theorem}
\newtheorem{lemma}{Lemma}
\newtheorem{proposition}{Proposition}

\theoremstyle{definition}
\newtheorem{definition}{Definition}

\theoremstyle{remark}
\newtheorem*{remark}{Remark}

\begin{document}

\title{Osserman Conjecture in dimension $n \ne 8, 16$}

\author{Y.Nikolayevsky}

\date{}

\maketitle

\begin{abstract}
Let $M^n$ be a Riemannian manifold and $R$ its curvature tensor. For a point
$p \in M^n$ and a unit vector $X \in T_pM^n$, the Jacobi operator is defined by $R_X =
R(X, \cdot \,) X$. The manifold $M^n$ is called {\it pointwise Osserman} if, for every
$p \in M^n$, the spectrum of the Jacobi operator does not depend of the choice of $X$,
and is called {\it globally Osserman} if it depends neither of $X$, nor of $p$.
Osserman conjectured that globally Osserman manifolds are two-point homogeneous. We
prove the Osserman Conjecture for $n \ne 8, 16$, and its pointwise version for
$n \ne 2,4,8,16$. Partial result in the case $n = 16$ is also given.
\end{abstract}

\section{Introduction}
\label{intro}

An {\it algebraic curvature tensor} $R$ in a Eucliean space $\Rn$ is a $(3, 1)$
tensor having the same symmetries as the curvature tensor of a Riemannian manifold.
For $X \in \Rn$, the {\it Jacobi operator} $R_X : \Rn \to \Rn$ is defined by
$R_XY = R(X, Y)X$ . The Jacobi operator is symmetric and $R_XX = 0$ for all
$X \in \Rn$. Throughout the paper, "eigenvalues of the Jacobi operator"
refers to eigenvalues of the restriction of $R_X$, with $X$ a
unit vector, to the subspace $X^\perp$.

\begin{definition} An algebraic curvature tensor $R$ is called
{\it Osserman} if the eigenvalues of the Jacobi operator $R_X$ do not
depend of the choice of a unit vector $X \in \Rn$.
\end{definition}

\begin{definition} A Riemannian manifold $M^n$ is called
{\it pointwise Osserman} if its curvature tensor is Osserman. If, in
addition, the eigenvalues of the Jacobi operator are constant on $M^n$,
the manifold $M^n$ is called {\it globally Osserman}.
\end{definition}

Two-point homogeneous spaces ($\Rn, \mathbb RP^n, S^n, H^n, \mathbb CP^n,
\mathbb CH^n, \mathbb HP^n, \mathbb HH^n, \linebreak[4]
\mathbb CayP^2$, and $\mathbb CayH^2$)
are globally Osserman, since the isometry group of each of this
spaces is transitive on its unit sphere bundle. Osserman \cite{13}
conjectured that the converse is also true:

\medskip
\noindent
{\bf Osserman Conjecture.} {\it A globally Osserman manifold is two-point
homogeneous.}
\medskip

%\begin{conjecture}[Osserman] A globally Osserman manifold is two-point
%homogeneous.
%\end{conjecture}

For manifolds of dimension $n \ne 4 k, \,\, k \ge 2$ the Osserman Conjecture is
proved by Chi \cite{5}. Further progress was made in \cite{7,9,11,12}.
We refer to \cite{6} for results on Osserman Conjecture in
semi-Riemannian geometry. The characterization of $p$-Osserman manifolds
(the averaging of the Jacobi operator over any $p$-plane
has constant eigenvalues) was given by Gilkey in \cite{8}:
any $p$-Osserman Riemannian manifold with $2 \le p \le n-2$ has constant sectional
curvature.

Our main result is the following Theorem.
\begin{theorem} \label{th1}
A globally Osserman manifold of dimension
$n \ne 8, 16$ is two-point homogeneous.
A pointwise Osserman manifold of dimension
$n \ne 2, 4, 8, 16$ is two-point homogeneous.
\end{theorem}

Note that in
dimension two, any Riemannian manifold is pointwise Osserman, but globally
Osserman manifolds are the ones having constant Gauss curvature.

In dimension
four, the Osserman Conjecture is proved in \cite{5}. However, there exist
pointwise Osserman four-dimensional manifolds that are {\it not} two-point
homogeneous (see \cite[Corollary 2.7]{9}).

A (pointwise or globally) Osserman manifold of dimension eight is known
to be two-point homogeneous in each of the following cases:
(i) the Jacobi operator has an eigenvalue of multiplicity at least $5$
\cite[Theorem 7.1]{9}, \cite[Theorem 1.2]{11};
(ii) the Jacobi operator has no more than two distinct eigenvalues
\cite[Theorem 2]{12}.

In dimensions sixteen, we have the following Theorem.

\begin{theorem} \label{th2} A (pointwise or globally) Osserman manifold $M^{16}$ is
two-point \nolinebreak[4] ho\-mogeneous if the Jacobi operator has no eigenvalues of 
multiplicity
$7, 8$ and $9$.
\end{theorem}

The paper is organized as follows. In Section~\ref{sec:2}, we consider algebraic curvature
tensors with Clifford structure. All of them have the Osserman property, and, in
the most cases, the converse is also true. This is the key statement of the paper
(Proposition~\ref{p1}). Moreover, in the cases covered by the Theorems, the existence of
the Clifford structure on a manifold implies that the manifold is two-point
homogeneous (Proposition~\ref{p2}). Further in Section~\ref{sec:2}, we give the proof of
the
both Theorems assuming Proposition~\ref{p1}. Section~\ref{sec:3} contains the proof of
Proposition~\ref{p1} modulo Propositions~\ref{p3}, \ref{p4} which are proved in
Sections~\ref{sec:4} and \ref{sec:5}, respectively.

\section{Manifolds with Clifford structure. Proof of the Theorems}
\label{sec:2}

Osserman algebraic curvature tensors with Clifford structure were introduced
by Gilkey \cite{7}, Gilkey, Swann, Vanhecke \cite{9}:

\begin{definition} \label{d3}
An algebraic curvature tensor $R$ in $\Rn$ has
a {\rm $\Cliff(\nu)$-structure} if
\begin{multline}\label{e1}
R(X, Y) Z = \la_0 (\< X, Z \> Y - \< Y, Z \> X) \\
+ \sum_{s=1}^\nu \tfrac 13 (\mu_s - \la_0) (2 \< J_sX, Y \> J_sZ +
\< J_sZ, Y \> J_sX - \< J_sZ, X \> J_sY),
\end{multline}
where $J_1, \dots , J_\nu$ are skew-symmetric orthogonal operators
satisfying the Hurwitz relations
$J_s J_q + J_q J_s = -2 \K_{qs} I_n$ and
$\mu_s \ne \la_0$.

A Riemannian manifold $M^n$ has a {\it $\Cliff(\nu)$-structure} if its
curvature tensor does.
\end{definition}

For skew-symmetric operators $J_1, \dots , J_\nu$ the Hurwitz
relations are equivalent to the fact that
$\< J_s X, J_q X \> = \K_{sq} \|X\|^2$ for all $X \in \Rn$.
Note that some of the $\mu_s$'s in \eqref{e1} can be equal.

The Jacobi operator of the algebraic curvature tensor $R$ with the Clifford
structure given by \eqref{e1} has the form
\begin{equation} \label{e2}
R_XY =  \la_0 (\| X\|^2 Y - \< Y, X \> X) +
\sum_{s=1}^\nu (\mu_s -\la_0) \< J_sX, Y \> J_sX,
\end{equation}
and the tensor $R$ can be reconstructed from \eqref{e2} using polarization
and the first Bianchi identity.

For any unit vector $X$, the Jacobi operator $R_X$ given by \eqref{e2}
has constant eigenvalues $\la_0, \la_1, \dots , \la_p$, where $\la_1, \dots ,
\la_p$ are the $\mu_s$'s without repetitions. The eigenspace corresponding to
the eigenvalue $\la_\a, \,\, \a \ne 0 $ is
$E_{\la_\a}(X) = \Span_{s: \mu_s = \la_\a} \!\! (J_sX)$, and the $\la_0$-eigenspace
is $E_{\la_0}(X)\!\! = \!\! (\Span (X, J_1X, \dots , J_\nu X))^\perp\!,$ provided $\nu < 
n-1$.
Hence a $\Cliff(\nu)$ algebraic curvature tensor (manifold) is Osserman
(pointwise Osserman, respectively).

We will show that, at the most cases, the converse is also true. Note however,
that there exists at least one Osserman algebraic curvature tensor having no
Clifford structure, namely the curvature tensor of the Cayley projective plane
and, up to a sign, of its hyperbolic dual (see the Remark at the end of
Section~\ref{sec:3}).

For a unit vector $X$, the vectors $J_1 X, \dots , J_\nu X$ are linearly
independent (even orthonormal) and are tangent to the unit sphere
$S^{n-1} \subset \Rn$ at $X$. From the Adams Theorem \cite{1} it follows
that $\nu \le \rho(n) - 1$, where $\rho(n)$ is the Radon number, defined as
follows: for $n= 2^{4a+b} c$ with $c$ odd integer and $0 \le b \le 3, \quad
\rho(n) = 2^b + 8 a$. Moreover, for every $\nu \le \rho(n) - 1$, there exist
the operators $J_1, \dots, J_\nu$ with the required properties
\cite{3,10}, and so there exist algebraic curvature tensors in $\Rn$
having a $\Cliff(\nu)$-structure.

In \cite{9} the following approach to the Osserman Conjecture was suggested:
\begin{itemize}
\item[(i)] show that Osserman algebraic curvature tensors have Clifford structure;
\item[(ii)] classify Riemannian manifolds having curvature tensor as in
\textup{(i)}.
\end{itemize}

Following this scheme, we derive Theorem~\ref{th1} and Theorem~\ref{th2}
from two Propositions below. We show that by topological reasons, the Jacobi
operator of an
Osserman algebraic curvature tensor must always have an eigenvalue of
multiplicity $m \ge n - \rho(n)$. Denote the sum of multiplicities of all
the other eigenvalues by $\nu = n - 1 - m \le \rho(n) -1$
(we use the same notation $\nu$ as in Definition~\ref{d3} aiming to find a
$\Cliff(\nu)$-structure for $R$).
Our proof works when the number $\nu$ is small enough compared to $n$.
When $n \ne 8, 16$ this is guaranteed by the fact that $\rho(n)$ is small,
but for $n = 16$ we need to impose extra conditions on the spectrum of the
Jacobi operator, as in Theorem~\ref{th2}.

\begin{proposition}
\label{p1}
Let $R$ be an Osserman algebraic curvature tensor in
$\Rn$. Let $m$ be the maximal multiplicity of the eigenvalues of its Jacobi
operator and $\nu = n - 1 - m$. If
$$
n \ge 3 \nu \quad \text{and} \quad n > \frac {(\nu + 1)^2}{4},
$$
then $R$ has a $\Cliff(\nu)$-structure.
\end{proposition}

\begin{proposition}[\cite{11}, Theorem 1.2]
\label{p2}
A Riemannian manifold $M^n$ with a \linebreak[4]
$\Cliff(\nu)$-structure is two-point homogeneous, provided that
\begin{itemize}
\item[(a)] $n \ne 2, 4, 8, 16$, or
\item[(b)] $n = 8, \, \nu < 3$, or
\item[(c)] $n = 16, \, \nu \ne 8$.
\end{itemize}
\end{proposition}

\begin{proof} [Proof of Theorem~\ref{th1} and Theorem~\ref{th2}] Both Theorems
follow from Propositions~\ref{p1} and \ref{p2} directly, if we can show that
the number $\nu$ defined in Proposition~\ref{p1} satisfies the inequalities
$n \ge 3 \nu, \,\, n > \tfrac {(\nu + 1)^2}{4}$.

These inequalities follow from the fact that $\nu \le \rho(n) - 1$. Indeed, as the
formula for $\rho(n)$ shows, for all
$n \ne 2, 4, 8, 16, \quad n \ge 3 (\rho(n) - 1), \,\, n > \tfrac {\rho(n)^2}{4}$.
For $n = 16$, we get $\nu \le 8$. The hypothesis of Theorem~\ref{th2} then implies that
$\nu \le 5$, and so $3 \nu \le 16, \,\, \tfrac {(\nu + 1)^2}{4} < 16$.

Hence it remains to show that $\nu \le \rho(n) - 1$. Let $M^n$ be a pointwise
Osserman manifold. Locally, in a neighbourhood of a generic point $x \in M^n$,
the Jacobi operator has a constant number of eigenvalues, with constant
multiplicities. Let the Jacobi operator have $p+1$ distinct eigenvalues,
with multiplicities $m_0, m_1, \dots, m_p$, respectively,
$m_0 + m_1 + \dots + m_p = n - 1$. Let $m = m_0$ be the maximal multiplicity
and $\nu = n - 1 - m$.

For a unit vector $X \in T_xM^n$, the eigenspaces of the Jacobi operator are
mutually orthogonal subspaces of $T_XS^{n-1}$, of dimension $m_0, m_1,
\dots, m_p$. This gives $p+1$ continuous plane fields in
the tangent bundle $TS^{n-1}$ of the unit sphere $S^{n-1} \subset T_xM^n$.

We follow the arguments of \cite[p. 216]{15}. Let
$f: S^{n-2} \to \SO(n-1)$ be the clutching map for the tangent bundle
$TS^{n-1}$, and $[f] \in \pi_{n-2} \SO(n-1)$ its homotopy class. The bundle
$TS^{n-1}$ admits $p+1$ continuous orthogonal plane fields of dimension $m_0, m_1,
\dots, m_p$, iff $[f]$ lies in the subgroup
$\pi_{n-2} (\SO(m_0) \times \SO(m_1) \times \dots \times \SO(m_p))$ of $\pi_{n-2} \SO(n-
1)$
defined by the inclusion map
$i: \SO(m_0) \times \SO(m_1) \times \dots \times \SO(m_p) \to \SO(n-1)$. Similarly, the
bundle $TS^{n-1}$ admits $\nu = n - m_0 - 1$ continuous orthonormal vector
fields, iff $[f]$ lies in the subgroup $\pi_{n-2} \SO(m_0)$ of
$\pi_{n-2} \SO(n-1)$ defined by the inclusion map
$i': \SO(m_0) \to \SO(n-1)$.

But the image of $i_*$ lies in the image of $i'_*$ since $\SO(m_0)$ is
the largest of the $\SO(m_\a)$, and so every $\SO(m_\a)$ can be homotoped in
$\SO(n-1)$ to lie inside $\SO(m_0)$.

Hence there exist $\nu$ vector fields on $S^{n-1}$, and the Adams Theorem
\cite{1} gives $\nu \le \rho(n) - 1$.
%\qedsymbol
\end{proof}

The remaining part of the paper is devoted to the proof of Proposition~\ref{p1}.

\section{Proof of Proposition~\ref{p1}}
\label{sec:3}

Let $\tilde R$ be an Osserman algebraic curvature tensor in $\Rn$ such that
the corresponding Jacobi operator has $p+1$ distinct eigenvalues
$\tilde \la_0, \tilde \la_1, \dots, \tilde \la_p$ with multiplicities
$m_0, m_1, \dots, m_p$, respectively, $m_0 + m_1 + \dots + m_p = n - 1$.
Let $m_0$ be the maximal multiplicity and $\nu = n - 1 - m_0$, the sum of
all the other multiplicities.

Consider an algebraic curvature tensor
$R = \tilde R - \tilde \la_0 R^1$, where $R^1$ is the curvature tensor
of the unit sphere. Then $R$ is still Osserman, with the Jacobi operator
having eigenvalues $\la_\a = \tilde\la_\a-\tilde\la_0$ with multiplicities
$m_\a$ for $\a =1, \dots, p$, respectively, and the eigenvalue $0$ with
multiplicity $m_0$. To prove Proposition~\ref{p1} it is sufficient to
show that $R$ has a $\Cliff(\nu)$-structure.

Let $\mu_1, \dots , \mu_\nu$ be the $\la_\a$'s counting the multiplicities,
that is, $\mu_1 = \dots = \mu_{m_1} = \la_1, \,
\mu_{m_1+1} = \dots = \mu_{m_1+m_2} = \la_2, \dots,
\mu_{\nu-m_p+1} = \dots = \mu_\nu = \la_p$.
In a Euclidean space $\rnu$, choose an orthonormal basis $e_1, \dots e_\nu$ and
define the operator $\L: \rnu \to \rnu$ by $\L \, e_s = \mu_s e_s, \, s = 1, \dots , \nu$.
The matrix of $\L$ is then $\diag\{\mu_1, \dots , \mu_\nu\}$.

Proposition~\ref{p1} follows from the two Propositions below.

\begin{proposition}
\label{p3}
Let $R$ be an Osserman algebraic curvature tensor in
$\Rn$. Let $0$ be the eigenvalue of its Jacobi operator with the maximal
multiplicity $m$, and $\nu = n - 1 - m$.

Assume that $n \ge 3 \nu$\frenchspacing. Then there exists a linear
map $M \! :\Rn \to \Hom (\rnu, \Rn), \linebreak[2]  X \to M_X$ such that the Jacobi
operator admits the following linear decomposition:
\begin{equation} \label{e3}
R_X = M_X \, \L \, M_X^t.
\end{equation}
The map $X \to M_X$ is determined uniquely up to a precomposition $X \to M_X N$ with
an element $N$ from the group $\O_{\L} = \{N : N \L \, N^t = \L \}$.
\end{proposition}

\begin{proposition}
\label{p4}
Let $R$ be an Osserman algebraic curvature tensor in
$\Rn$ with the Jacobi operator having the form \eqref{e3}.
Assume that $n > \tfrac {(\nu + 1)^2}{4}$. Then $R$ has a $\Cliff(\nu)$-structure.
\end{proposition}

Proposition~\ref{p3} is proved in Section~\ref{sec:4}. Using the Osserman property
we successively show
that for $k = 1, 2, 3, n$ the following holds: for almost any set of $k$ orthonormal
vectors $E_1, \dots , E_k$ in $\Rn$, there exist linear operators
$M_1, \dots, M_k: \rnu \to \Rn$ such that $R_{x_1E_1 + \dots + x_kE_k}
= (x_1M_1 + \dots + x_kM_k) \, \L \, (x_1M_1 + \dots + x_kM_k)^t$ for all
$x_1, \dots , x_k$.
Then for a vector $X = x_1E_1 + \dots + x_nE_n$ we define
$M_X = x_1M_1 + \dots + x_nM_n$.

Proposition~\ref{p4} is proved in Section~\ref{sec:5}. We show that, with an
appropriate choice of
the basis $e_1, \dots e_\nu$ in $\rnu$, the operators $J_s$ in $\Rn$ defined
by $J_sX = M_X e_s$ give the Clifford structure for $R$.

\begin{remark} The claim of Proposition~\ref{p1} fails to be true at least in the
case when $n = 16, \,\, \nu = 7$, since
the curvature tensor of the Cayley projective plane $\mathbb Cay P^2$
(and of its hyperbolic dual $\mathbb Cay H^2$) has no Clifford structure.

This follows from the fact that, unlike the holonomy groups of
$\mathbb CP^n$ and $\mathbb HP^n$, the holonomy group
$\operatorname {Spin}(9)$ of the Cayley projective plane
has no proper normal subgroups \cite{2,4}.
The nonexistence of the Clifford structure is also confirmed by
the following octonionic computation based on the formula for the
curvature tensor of $\mathbb Cay P^2$ \cite[Theorem 6.1]{4}.

Identify a tangent space to $\mathbb Cay P^2$ with $\mathbb Cay \oplus \mathbb Cay$.
Then for orthogonal vectors $X = (a, b), \,\, Y = (c, d)$
the Jacobi operator has the form
$$
R_XY = \frac {\a}{4} ((4 \|a\|^2 + \|b\|^2) c + 3 (ab)d^*,
(4 \|b\|^2 + \|a\|^2) d + 3 c^* (a b)),
$$
where ${}^*$ is the octonion conjugation and $\|a\|^2 = a a^*,\,\,
\<a, b\> = \tfrac12 (a b^* + b a^*)$.

It follows that for any unit vector $X$ the Jacobi operator $R_X$
has two eigenvalues: $\a$, of multiplicity $7$, with the eigenspace
$$
E_\a(X) = \{(c, d) : ad = cb, \,\, \< a, c\> = \<b, d\> = 0\},
$$
and $\tfrac {\a}{4}$, of multiplicity $8$, with the eigenspace
$$
E_{\frac {\a}{4}}(X) = \{(c, d) : a (\|b\|^2 d - \<b, d\> b) =
(\|a\|^2 c - \<a, c\> a) b\}.
$$

As it follows from \eqref{e2}, the existence of a Clifford structure would
imply the existence of seven (respectively, eight) linear operators $J_i$ such
that $E_\a(X) = \Span (J_1X, \dots, J_7X)$
(respectively, $E_{\frac {\a}{4}}(X) = \Span (J_1X, \dots, J_8X)$).

However, it is not difficult to see that there is no nonzero
$\mathbb R$-linear operator $J: \mathbb Cay \oplus \mathbb Cay \to \mathbb Cay \oplus
\mathbb Cay$
such that for all $X, \quad JX \in E_\a(X)$. With some calculation, one can show
that the same is true for $E_{\frac {\a}{4}}(X)$, as well.

Thus the curvature tensor of $\mathbb Cay P^2$ admits neither $\Cliff(7)$-, nor
$\Cliff(8)$-
structure, hence no Clifford structure at all.
\end{remark}

\section{Proof of Proposition~\ref{p3}}
\label{sec:4}

The proof goes by the following plan. First, in Lemma~\ref{l1}, we show that for any
vector $X$,
there exists an operator $M_X$ satisfying \eqref{e3}. Next, in
Lemma~\ref{l3}, we prove that for almost any two vectors $X, Y$, the operators $M_X$
and $M_Y$ can be chosen accordingly, that is, in such a way that
$R_{x X + y Y} = (M_X x + M_Y y) \, \L \, (M_X x + M_Y y)^t$ for all real $x, y$.
In Lemma~\ref{l4} we extend this result to the case of three vectors.
Then we show that the existence of the linear decomposition
of the form \eqref{e3} for almost any three vectors already implies the
existence of the global decomposition.

\medskip

Let $R$ be an Osserman algebraic curvature tensor in $\Rn$,
with the Jacobi operator having $p$ distinct nonzero eigenvalues
$\la_1, \la_2, \dots, \la_p$ of
multiplicities $m_1, m_2, \dots,   m_p$, respectively, and the eigenvalue
$0$ of multiplicity $n - 1 - \nu$.

For a nonzero vector $X$, the subspaces $\Ker R_X$ and $\Im R_X$ are orthogonal
and have dimension $n - \nu$ and $\nu$, respectively. By the Raki\`c duality principle
\cite{14}, for two orthonormal vectors $X, Y$, the following holds: $Y$ is an
eigenvector of $R_X$, if and only if $X$ is an eigenvector of $R_Y$ (with the
same eigenvalue). For the eigenvalue $0$, we get
\begin{equation}\label{e4}
Y \in \Ker R_X \quad \text{iff} \quad X \in \Ker R_Y.
\end{equation}
We will use a slight modification of the duality principle, noting that for
the eigenvalue $0$ the assumption of orthonormality of $X$ and $Y$ can be dropped.
Indeed, let $\psi \ne 0, \tfrac {\pi}2, \pi$ be the angle between unit vectors $X, Y$, and
let $Z$ be a unit vector in $\Span(X, Y)$ orthogonal to $X$ and such that
$Y = \cos \psi X + \sin \psi Z$. Assume that $Y \in \Ker R_X$. Since
$X \in \Ker R_X$, we have $Z \in \Ker R_X$, and so $X \in \Ker R_Z$. Then
$R_YX = R_{\cos \psi X + \sin \psi Z}X =
\cos \psi \sin \psi \,\,R(Z,X)X + \sin ^2 \psi \,\, R_ZX = 0$, that is,
$X \in \Ker R_Y$.

\medskip

We first show that the decomposition claimed in the Proposition exists for
every single operator $R_X$.

\begin{lemma} \label{l1}
For any unit vector $X$, there exists a linear operator $M_X: \rnu \to \Rn$ such that
$$
R_X = M_X \, \L \, M_X^t.
$$
The operator $M_X$ is determined uniquely up to a precomposition with an
element $N \in \O_{\L} = \{N : N \L \, N^t = \L \}$. Moreover,
$\Im R_X = \Im M_X$ and $M_X^tM_X = N^tN$ for some $N \in \O_{\L}$.
\end{lemma}

\begin{proof} Let $E_1, \dots , E_\nu$ be an orthonormal basis of
eigenvectors of $R_X$ with nonzero eigenvalues. For a vector
$y = (y_1, \dots, y_\nu) \in \rnu$ define
$M_X y = y_1 E_1 + \dots + y_\nu E_\nu$. The symmetric
operators $R_X$ and $M_X \L M_X^t$
acting in $\Rn$ have the same eigenvectors and the same eigenvalues, hence
$R_X = M_X \L M_X^t$. Moreover, $M_X^tM_X = I_{\nu}$, since for any
$y, z \in \rnu$ we have $\<M_X y, M_X z \> = \<y, z\>$.

If $\tilde M_X: \rnu \to \Rn$ is another operator such that
$R_X = \tilde M_X \L \tilde M_X^t$, then $\tilde M_X \L \tilde M_X^t = M_X \L M_X^t$,
and so $\tilde M_X = M_X \L M_X^t \tilde M_X (\tilde M_X^t \tilde M_X)^{-1} \L^{-1}$
(the operator $\tilde M_X^t \tilde M_X$ is nonsingular since
$\rk \tilde M_X = \rk R_X = \nu$). Therefore $\tilde M_X = M_X N$, for some
operator $N$ in $\rnu$ and $M_X (\L - N \L N^t) M_X^t = 0$. Since
$\rk M_X = \nu, \linebreak[3] N \in \O_{\L}$.
%\qedsymbol
\end{proof}

Next we need the following generic position Lemma. Denote $V_k (\Rn)$ the Stiefel
manifold of $k$-tuples of orthonormal vectors in $\Rn$.

\begin{lemma} \label{l2}
1. Let $n \ge 2 \nu$. Then the set
$S_2 = \{(X, Y) \in V_2 (\Rn) \, : \, \Im R_X \cap \Im R_Y = 0\}$
is open and dense in $V_2 (\Rn)$.

2. Let $n \ge 3 \nu$. Then the set
$S_3 = \{(X, Y, Z) \in V_3 (\Rn) \, : \,
\dim (\Im R_X + \Im R_Y + \Im R_Z) = 3 \nu\}$ is open and dense in $V_3(\Rn)$.

3. If $(X, Y, Z) \in S_3$, then for any unit vector $U \in \Span (Y, Z)$
the pair $(X, U)$ is in $S_2$.
\end{lemma}

\begin{proof} The proof of 1. and 2. is quite similar and is based
on the dimension count. Both $S_2$ and $S_3$ are open.
We claim that they are also dense.

1. Let $X$ be a unit vector in $\Rn$, and $S^{n-2}$ be the unit sphere in the
subspace $X^\perp$. We want to show that for an open dense set of
vectors $Y \in  S^{n-2}, \quad \dim(\Ker R_X \cap \Ker R_Y) = n - 2 \nu$.

Let $S$ be the unit sphere in the subspace $\Ker R_X, \quad
\dim S = n-\nu-1$, and let $E$ be a vector bundle with the base $S$
and the fiber $F_Z = \Ker R_Z \cap X^\perp$ over a point
$Z \in S$ ($\dim F_Z = n - \nu - 1$ since $X \in \Ker R_Z$ by the duality
principle \eqref{e4}). Then $SE$, the corresponding unit sphere bundle, is a
compact analytic manifold of dimension $2 n - 2 \nu -3$.

Define the projection map $\pi: SE \to S^{n-2}$ by $\pi(Z, Y) = Y$.
By the duality principle, for any
$Y \in S^{n-2}, \quad \pi^{-1}(Y) =\{(Z, Y) \in SE : Z \in \Ker R_Y\}$.
The map $\pi$ is differentiable (even analytic) since the
subspace $\Ker R_Z$ viewed as a point of the corresponding Grassmannian
depends analytically on $Z$.

Now if $n = 2 \nu$, then $2 n - 2 \nu - 3 < n - 2$, and so for every $Y$ from
the open dense subset $S^{n-2} \setminus \pi(SE)$ of $S^{n-2}, \quad
\Ker R_Y \cap \Ker R_X = 0$.

If $n > 2 \nu$, then $\Ker R_Y \cap \Ker R_X \ne 0$ for any $Y$, hence the
map $\pi$ is surjective. By the Sard Theorem, for an open dense set of the
$Y$'s in $S^{n-2}, \quad d\pi$ has the maximal rank $n - 2$ at all the points of
$\pi^{-1}(Y)$. For such points $\pi^{-1}(Y) = (S_Y^{n - 2\nu - 1}, Y)$,
where $S_Y^{n - 2\nu - 1}$ is the unit sphere in $\Ker R_X \cap \Ker R_Y$. So
$\dim(\Ker R_X \cap \Ker R_Y) = n - 2 \nu$.

2. Let $(X, Y) \in S_2$, and let $S^{n-3}$ be the unit sphere in the subspace
$X^\perp \cap Y^\perp$. We show that for an open dense set of vectors
$Z \in S^{n-3}, \, \dim(\Ker R_X \cap \Ker R_Y \cap \Ker R_Z) = n - 3 \nu$.

Let $S$ be the unit sphere in the subspace $\Ker R_X \cap \Ker R_Y,
\,\, \dim S = n-2\nu-1$, and let $E$ be a vector bundle with the base $S$
and the fiber $F_U = \Ker R_U \cap (X^\perp \cap Y^\perp)$ over a point
$U \in S$ ($\dim F_U = n - \nu - 2$ since $X, Y \in \Ker R_U$). Then
$SE$, the corresponding unit sphere bundle, is a compact analytic manifold of
dimension $2 n - 3 \nu -4$.

The projection map $\pi: SE \to S^{n-3}$ defined by $\pi(U, Z) = Z$ is
also analytic. By the duality principle, for any $Z \in S^{n-3}, \quad
\pi^{-1}(Z) =\{(U, Z) \in SE : U \in \Ker R_Z \}$.

If $n = 3 \nu$, then the image of $\pi$ does not cover $S^{n-3}$, and we can
take any $Z$ from its complement.

Otherwise, $\pi$ is surjective. Applying the Sard Theorem, we find an open dense
set of points $Z \in S^{n-3}$, such that $\pi^{-1}(Z) = (S_Z^{n - 3\nu-1}, Z)$,
where $S_Z^{n - 3\nu-1}$ is the unit sphere in
$\Ker R_X \cap \Ker R_Y \cap \Ker R_Z$. So
$\dim(\Ker R_X \cap \Ker R_Y \cap \Ker R_Z) = n - 3 \nu$.

3. By the duality principle \eqref{e4},
$\Ker R_U \supset (\Ker R_Y \cap \Ker R_Z)$, and so
$\Im R_U \subset (\Im R_Y + \Im R_Z)$. Then
$\Im R_U \cap \Im R_X \subset (\Im R_Y + \Im R_Z) \cap \Im R_X = 0$, since
$(X, Y, Z) \in S_3$.
%\qedsymbol
\end{proof}

We now show that the operator $R_X$ admits the linear decomposition \eqref{e3}
on almost every two-plane in $\Rn$. Define the symmetric operator $R_{XY}$ by
$R_{XY}Z = \frac 12 (R(X, Z) Y + R(Y,Z) X)$.

\begin{lemma} \label{l3}
Suppose that $n \ge 2 \nu$. Then for any pair of orthonormal
vectors $(X, Y) \in S_2$, there exist linear operators $M_1, M_2: \rnu \to \Rn$
such that  for all $x, y \in \mathbb R$,
$$
R_{x X + y Y} = (M_1 x + M_2 y) \, \L \, (M_1 x + M_2 y)^t.
$$
Moreover, the operators $M_1, M_2$ are determined  uniquely up to a
precomposition $M_1 N, \, M_2 N$ with an element $N \in \O_{\L}$.
\end{lemma}

The uniqueness part can be rephrased as follows: once $M_1$ with the property
$M_1 \L M_1^t = R_X$ is chosen, then there exists a unique $M_2$ such
that the pair $M_1, M_2$ satisfies the equation of Lemma~\ref{l3}.

\begin{proof} We have $R_{xX + yY} = R_X x^2 + C x y + R_Y y^2$, with
a symmetric operator $C = 2 R_{XY}$. The claim is equivalent to the fact that
$$
R_X = M_1 \L M_1^t, \quad R_Y = M_2 \L M_2^t, \quad
C = M_1 \L M_2^t + M_2 \L M_1^t.
$$
By Lemma~\ref{l1} we can find two operators, $M_1, M_2: \rnu \to \Rn$, such that
$R_X = M_1 \L M_1^t, \quad R_Y = M_2 \L M_2^t$, and so our goal is to show that
they can be chosen in such a way that $C = M_1 \L M_2^t + M_2 \L M_1^t$.

By the Osserman property, the operator $R_{\cos \phi X + \sin \phi Y}$ is
isospectral, for all $\phi \in \mathbb R$. Its eigenspaces (viewed as the curves
in the corresponding Grassmannians) are analytic with respect to $\phi$.
Locally, in a neighbourhood of the point $\phi = 0$, there exists an analytic
orthogonal transformation $U(\phi)$ such that
$U(0) = I_n$ and $R_{\cos \phi X + \sin \phi Y} = U (\phi) R_X U (\phi)^t$.
Let $U(\phi) = I_n + K \phi + (\frac 12 K^2 + K_1) \phi^2 + o(\phi^3)$ be
the Taylor expansion at $\phi = 0$, with
$K$ and $K_1$ skew-symmetric operators.

Then we have
\begin{equation}\label{e5}
C = [K, R_X], \quad R_Y = R_X + \tfrac 12 [K, C] + [K_1, R_X],
\end{equation}
and so $C = \tilde M_2 \L M_1^t + M_1 \L \tilde M_2^t$, with $\tilde M_2 = K M_1$.

Similar arguments applied at the point $\phi = \pi/2$ show that
$C = \tilde M_1 \L M_2^t + M_2 \L \tilde M_1^t$ for some operator
$\tilde M_1: \rnu \to \Rn$. Equating the expressions for $C$ we get
\begin{equation}\label{e6}
\tilde M_1 \L M_2^t + M_2 \L \tilde M_1^t
= \tilde M_2 \L M_1^t + M_1 \L \tilde M_2^t.
\end{equation}
Let $Z \in \Ker R_X \cap \Ker R_Y$. Then $M_1^tZ = M_2^tZ = 0$. Acting on the
vector $Z$ by the both sides of \eqref{e6} we get
$M_1 (\L \tilde M_2^tZ) = M_2 (\L \tilde M_1^tZ)$. The pair $(X, Y)$ was chosen
in $S_2$, so the subspaces $\Im M_1 = \Im R_X$ and
$\Im M_2 = \Im R_Y$ have zero intersection in $\Rn$. It follows that
$\tilde M_2^tZ = \tilde M_1^tZ = 0$. Hence
$\Im \tilde M_1, \Im \tilde M_2 \subset \Im M_1 \oplus \Im M_2$.
In other words, there
exist linear operators $S_1, S_2, S_3, S_4 : \rnu \to \rnu$ such that
$$
\tilde M_1 = M_1 S_1 + M_2 S_3, \quad
\tilde M_2 = M_1 S_2 + M_2 S_4.
$$
Substituting this back to \eqref{e6} we get
\begin{multline*}
M_2 (S_3 \L + \L S_3^t)M_2^t + M_2 (\L S_1^t- S_4 \L) M_1^t \\
- M_1 (S_2 \L + \L S_2^t) M_1^t + M_1 (S_1  \L - \L S_4^t)M_2^t  = 0.
\end{multline*}
Using again the fact that $\Im M_1 \cap \Im M_2 = 0$ we find that $S_2 \L, S_3 \L$ are
skew-symmetric operators in $\rnu$, and $\L S_1^t = S_4 \L$. Then
$$
C = M_2 S_4 \L M_1^t + M_1 \L S_4^t M_2^t.
$$
Take a vector $Z$ in $\Ker R_X$. Then $M_1^t Z = 0$, and
$M_1^t K Z = - (S_2^t M_1^t + S_4^t M_2^t)Z = - S_4^t M_2^t Z$
since $\tilde M_2 = K M_1 = M_1 S_2 + M_2 S_4$.
Acting on $Z$ by the both sides of the second equation of \eqref{e5}, and
then taking the inner product with $Z$, we obtain $\<R_YZ, Z \> = \tfrac 12 \<[K, C]Z, Z\>
= -\<C K Z, Z\>$ since $R_XZ = 0$. Substituting $R_Y = M_2 \L M_2^t, \,\,
C = M_2 S_4 \L M_1^t + M_1 \L S_4^t M_2^t$ we find
\begin{multline*}
\<M_2 \L M_2^t Z, Z \> = - \<(M_2 S_4 \L M_1^t + M_1 \L S_4^t M_2^t) K Z, Z\> \\
 = - \<M_2 S_4 \L M_1^t K Z, Z\> - \< \L S_4^t M_2^t K Z, M_1^tZ\>
= \<M_2 S_4 \L S_4^t M_2^t Z, Z\>,
\end{multline*}
and so $\<(S_4 \L S_4^t - \L) (M_2^t Z), (M_2^t Z) \> = \<M_2 (S_4 \L S_4^t - \L)
M_2^t Z, Z \> = 0$
for all $Z \in \Ker R_X$.

The restriction of the operator $M_2^t$ to $\Ker R_X$ is epimorphic
(otherwise the images of $R_X$ and $R_Y$ would have a nonzero intersection),
so $M_2^t(\Ker R_X) = \rnu$. Hence the symmetric operator $S_4 \L S_4^t - \L$
vanishes. It follows that $S_4 \in \O_\L$.

Now replace $M_2$ by $\hat M_2 = M_2 S_4$. Then $\hat M_2 \L \hat
M_2^t$ is still $R_Y$, and \linebreak[4] $C = M_1 \L \hat M_2^t + \hat M_2 \L M_1^t$.
So $M_1, \, \hat M_2$ is the sought pair of operators.

To finish the proof it remains to show the uniqueness.
Suppose that for operators $M_1, M_2, M_3, M_4 : \rnu \to \Rn$,
$$
R_{x X + y Y} = (M_1 x + M_2 y) \, \L \, (M_1 x + M_2 y)^t =
(M_3 x + M_4 y) \, \L \, (M_3 x + M_4 y)^t.
$$
By Lemma~\ref{l1}, there exist $N_1, N_2 \in \O_\L$ such that
$M_3 = M_1 N_1, \, M_4 = M_2 N_2$. Equating the terms with $x y$ we then obtain
$$
M_1 (\L - N_1 \L N_2^t) M_2^t + M_2 (\L - N_2 \L N_1^t) M_1^t = 0.
$$
Since $\Im M_1 \cap \Im M_2 = 0$, it follows that $N_1 \L N_2^t = \L$, hence
$N_1 = N_2$.
%\qedsymbol
\end{proof}

The next step is to show that $R_X$ admits the linear decomposition \eqref{e3}
on almost every three-space in $\Rn$.

\begin{lemma}
\label{l4}
Suppose that $n \ge 3 \nu$. Then for any triple of orthonormal
vectors $(X, Y, Z) \in S_3$, there exist linear operators
$M_1, M_2, M_3: \rnu \to \Rn$ such that for all $x, y, z \in \mathbb R$,
$$
R_{x X + y Y + z Z} = (M_1 x + M_2 y + M_3 z) \, \L \, (M_1 x + M_2 y + M_3 z)^t.
$$
\end{lemma}

\begin{proof} Since the triple $(X, Y, Z)$ is in $S_3$, every pair
$(X, \cos \phi Y + \sin \phi Z)$ must be in $S_2$ by 3. of
Lemma~\ref{l2}. Then for any $\phi \in \mathbb R$, we can find the operators
$M_1, M_2(\phi)$ such that
\begin{equation}\label{e7}
R_{x X + y (\cos \phi Y + \sin \phi Z)} =
(M_1 x + M_2(\phi) y) \, \L \, (M_1 x + M_2(\phi) y)^t,
\end{equation}
(the fact that $M_1$ can be chosen independent of $\phi$ follows from Lemma~\ref{l1}).

Denote $M_2 = M_2(0), \, M_3 = M_2(\pi/2)$. Then
$$
R_X = M_1 \L M_1^t, \quad R_Y = M_2 \L M_2^t, \quad R_Z = M_3 \L M_3^t.
$$
Let $P(\phi) = M_2(\phi) - \cos \phi M_2 - \sin \phi M_3$. We want to show that
$P(\phi) = 0$.

The terms with $x y$ of \eqref{e7} give
$M_1 \L P(\phi) ^t + P(\phi) \L M_1^t = 0$, and so $P(\phi) = M_1 S(\phi)$, with
$S(\phi) = - \L P(\phi) ^t M_1 (\L M_1^tM_1)^{-1}$ linear operator in $\rnu$.
Substituting $P(\phi) = M_1 S(\phi)$ in the terms of
\eqref{e7} with $y^2$, and dividing by $\cos \phi \sin \phi$ we obtain
\begin{multline*}
M_1(S(\phi) \L M_2(\phi)^t (\cos \phi \sin \phi)^{-1}) +
M_2 (\L P(\phi)^t /\sin \phi) \\ +
M_3 (\L P(\phi)^t /\cos \phi) = 2 R_{YZ} - M_2 \L M_3^t - M_3 \L M_2^t.
\end{multline*}
The fact that $(X, Y, Z) \in S_3$ means that the images
$\Im M_1 = \Im R_X, \linebreak[2] \Im M_2 = \Im R_Y$, and $\Im M_3 = \Im R_Z$
span a subspace of dimension $3 \nu$ in $\Rn$. It follows that the
operator $M_1 \oplus M_2 \oplus M_3 : \rnu \to \Rn$ is one-to-one and so
operators $S(\phi) \L M_2(\phi)^t /(\cos \phi \sin \phi), \, \L P(\phi)^t /\sin \phi$
and $\L P(\phi)^t /\cos \phi$ are independent of $\phi$. In particular,
both $P(\phi)/\sin \phi$ and $P(\phi)/\cos \phi$ must be constant. This is only
possible when $P(\phi)$ vanishes identically.

Then $M_2(\phi) = \cos \phi \, M_2 + \sin \phi \, M_3$ and the claim follows
from \eqref{e7}.
%\qedsymbol
\end{proof}

With Lemma~\ref{l4}, we can finish the proof of the Proposition as follows.

Choose an orthonormal basis $E_1, \dots , E_n$ in $\Rn$ in such a way that every
triple $(E_i, E_j, E_k)$ is in $S_3$ and every pair $(E_i, E_j)$ is in
$S_2$.
The set of such bases is open and dense in the Stiefel manifold
$V_n(\Rn) = \O(n)$.

For every $i = 2, \dots, n$, let $M_1, \, M_i$ be the operators constructed
as in Lemma \ref{l3} on the vectors $E_1, E_i$ (by Lemma~\ref{l1} we can take $M_1$
the same for all the $i$'s). Then for $i = 2, \dots , n$
$$
M_1 \L M_1^t = R_1, \quad
M_i \L M_i^t = R_i, \quad
M_1 \L M_i^t + M_i \L M_1^t = 2 R_{1i},
$$
where $R_k = R_{E_k}, \,\, R_{kl} = R_{E_kE_l}$.

By Lemma~\ref{l4}, for any pair $i \ne j, \, i,j \ge 2$ there exist operators
$\tilde M_i, \tilde M_j$ satisfying
\begin{equation}\label{e8}
R_{x E_1 + y E_i + z E_j} = (M_1 x + \tilde M_i y + \tilde M_j z) \, \L \,
(M_1 x +  \tilde M_i y + \tilde M_j z)^t.
\end{equation}
In particular,
$\tilde M_i \L \tilde M_j^t + \tilde M_j \L \tilde M_i^t  = 2 R_{ij}$. On the
other
hand, taking $z = 0$ in \eqref{e8} and applying the uniqueness part of Lemma~\ref{l3}
we get $\tilde M_i = M_i$. Similarly, $\tilde M_j = M_j$. It follows that
$M_i \L M_j^t + M_j \L M_i^t = 2 R_{ij}$.

Now for an arbitrary vector $X = x_1 E_1 + \dots + x_n E_n$, define the operator
$M_X: \rnu \to \Rn$ by
$$
M_X = x_1 M_1 + \dots + x_n M_n.
$$

Then $M_X \L M_X^t = \sum_{i=1}^n R_i x_i^2 + \sum_{i<j} 2R_{ij} x_i x_j = R_X$.

The fact that the map $X \to M_X$ is determined uniquely up to a precomposition
with a fixed element from $\O_\L$ follows from Lemma~\ref{l2}.

\section{Proof of Proposition~\ref{p4}}
\label{sec:5}

We are given an Osserman algebraic curvature tensor $R$ in $\Rn$
with the Jacobi operator having $p+1$ distinct eigenvalues,
$\la_1, \la_2, \dots, \la_p$, and $0$, of multiplicities
$m_1, m_2, \dots, m_p$, and $n - 1- \nu$, respectively.
The number $\nu$ satisfies the inequality
\begin{equation}\label{e9}
n > \frac {(\nu + 1)^2}{4} .
\end{equation}

In the Euclidean space $\rnu$, with the fixed orthonormal basis
$e_1, \dots, e_\nu$, the linear operator $\L$ is defined by
$\L e_s = \mu_s e_s, \, s = 1, \dots, \nu$, where
$\mu_1 = \dots = \mu_{m_1} = \la_1, \, \mu_{m_1+1} = \dots = \mu_{m_1+m_2} =
\la_2, \,\,\dots,\,\, \mu_{\nu-m_p+1} = \dots = \mu_\nu = \la_p$.

The Jacobi operator of $R$ has the form
\begin{equation}\label{e10}
R_X = M_X \, \L \, M_X^t,
\end{equation}
where $M:\Rn \to \Hom (\rnu, \Rn), \; X \to M_X$ is a linear map determined
uniquely up to a precomposition with an element $N$ from the group
$\O_\L \! = \! \{N \! : \! N \L N^t \! = \! \L\}$.

The central role in the proof is played by a quadratic map
$\Phi: \Rn \to \Hom(\rnu, \rnu)$ defined by
$$
\Phi(X) = M_X^tM_X, \qquad X \in \Rn.
$$

In terms of the map $\Phi$, the Osserman property of $R$ has the following
form.

\begin{lemma} \label{l5}
For every unit vector $X$,

1. the operator $\L \, \Phi(X): \rnu \to \rnu$ is similar to $\L$; in particular,
it has the same spectrum as $\L$;

2. there exists $N \in \O_\L$ such that $N^t \Phi(X) N = I_\nu$.
\end{lemma}

\begin{proof}By Lemma~\ref{l1}, for every unit vector $X$, there exists an
element $N$ (depending on $X$) in the group $\O_\L$ such that $\Phi(X) = N^t N$.
Then $\L \, \Phi(X) = \L N^t N = N^{-1} \L N$.
%\qedsymbol
\end{proof}

\medskip

The proof of the Proposition goes by induction by $p$, the number of distinct
nonzero eigenvalues of the Jacobi operator.

{\bf Base.} Let $p = 1$, that is, the Jacobi operator has only two eigenvalues:
$\la_1$ with multiplicity $\nu$, and $0$ with multiplicity $n - 1 - \nu$.
Then $\L = \la_1 I_\nu$ and by 1. of Lemma~\ref{l5}
$$
M_X^t M_X = \Phi(X) = \|X\|^2 I_\nu.
$$

Define the operators $J_s: \Rn \to \Rn, \quad s = 1, \dots, \nu$ by
$J_s X = M_X e_s$. Then for all $X, Y \in \Rn$,
$$
R_XY = \la_1 M_X \, M_X^t Y = \la_1 \sum_{s=1}^\nu \< M_X^t Y, e_s \> M_X e_s
= \la_1 \sum_{s=1}^\nu \< J_s X, Y \> J_s X,
$$
that is, the Jacobi operator has the required form (\eqref{e2}, with $\la_0 =
0$).

Moreover, for any nonzero $X$, the vectors $J_1 X, \dots , J_\nu X$ are linearly
independent (otherwise $\rk R_X < \nu$), and so all the operators $J_s$ are
skew-symmetric since $R_XX=0$.

We also have $J_s J_q + J_q J_s = - 2 \K_{qs} I_n$ for all
$1 \le q, s \le \nu$, since for any vector
$X, \; \<J_s X, J_q X \> = \<M_X e_s, M_X e_q\> = \<M_X^t M_X e_s, e_q\>
 =\|X\|^2 \K_{qs}$.

Thus the Jacobi operator has the form \eqref{e2}, with the skew-symmetric
orthogonal operators $J_1, \dots , J_\nu$ satisfying the Hurwitz relations.

{\bf Step.} The plan of proof is the following. Suppose that we already know
that for $p = k \ge 1$ any Osserman algebraic curvature tensor has a Clifford
structure. Let $R$ be an Osserman curvature tensor with $p = k+1$ distinct
nonzero eigenvalues of the Jacobi operator.

For every unit vector $X$, the $\la_\a$-eigenspace of $R_X$ is
$E_{\la_\a}(X) = \{M_Xu : \linebreak[4]
u \in \rnu, \, R_X (M_Xu) = \la_\a (M_Xu)\}$. This
defines a subspace $S_{\la_\a}(X) \subset \rnu$ of dimension $m_\a$, the
multiplicity of the eigenvalue $\la_\a$, consisting of vectors $u \in \rnu$
satisfying $R_X (M_Xu) = \la_\a (M_Xu)$.

The key step is to show that, with a particular choice of $\la_\a$,
the subspace $S_{\la_\a}(X)$ is independent of $X$ (Lemma~\ref{l7}), that is, there
exists a fixed subspace $S \subset \rnu, \,\, \dim S = m_\a$ such that
$E_{\la_\a}(X) = M_XS$ for all unit vectors $X \in \Rn$.
We then choose a basis $u_1, \dots, u_{m_\a}$ in $S$ and define the
operators $J_s: \Rn \to \Rn, \quad s = 1, \dots, m_\a$ by
$$
J_s X = M_X u_s.
$$
For every unit vector $X$ and for every $s, \quad J_sX$ is an eigenvector of
$R_X$ with the eigenvalue $\la_\a$. In particular, $\<J_sX, X\> = 0$ and so
all the $J_s$'s are skew-symmetric. In Lemma~\ref{l8} we show that the basis
$u_1, \dots, u_{m_\a}$ can be chosen in such a way that the operators $J_s$
are also orthogonal and satisfy the Hurwitz relations.

Introduce an algebraic curvature tensor $\hat R$ defined by its Jacobi operator
as
\begin{equation}\label{e11}
\hat R_XY = \la_\a \sum_{s=1}^{m_\a} \< J_s X, Y \> J_s X.
\end{equation}
Then $\hat R$ is Osserman, with the Jacobi operator having two eigenvalues,
$\la_\a$ and $0$.

Moreover, for every unit vector $X$, the $\la_\a$-eigenspace of $R_X$ and
$\hat R_X$ is the same: $\Span (J_1X, \dots , J_{m_\a}X)$. It follows that
the algebraic curvature tensor $R - \hat R$ is also Osserman. Its
Jacobi operator $R_X - \hat R_X$, for any unit vector $X$, has constant
eigenvalues $\la_1 , \dots, \la_{\a-1}, \la_{\a+1}, \dots, \la_{k+1}, 0$
with constant multiplicities (in fact, the $\la_\b$-eigenspaces of
$R_X - \hat R_X$ are the same as that of $R_X$, and
$\Ker (R_X - \hat R_X) = \Ker R_X \oplus E_{\la_\a}(X)$).

The number of nonzero eigenvalues of the Jacobi operator of $R - \hat R$
is one less than that for $R$, and so by the induction assumption the
algebraic curvature tensor $R - \hat R$ has a Clifford structure:
\begin{equation}\label{e12}
(R_X - \hat R_X)Y = \sum_{i=1}^{\nu-m_\a} \mu_i \< J_i X, Y \> J_i X
\end{equation}
for all $X, Y \in \Rn$, with skew-symmetric orthogonal operators $J_i$
satisfying the Hurwitz relations. Together with \eqref{e11} this gives a
$\Cliff(\nu)$-structure for $R$, provided the operators $J_s$ in
\eqref{e11} and the operators $J_i$ in \eqref{e12} satisfy the Hurwitz
relations $J_iJ_s + J_sJ_i = 0$. This is indeed the case, since for any
$X \in \Rn, \; J_sX$ and $J_iX$ are eigenvectors of $R_X$ which
correspond to different eigenvalues, and so are orthogonal.

This proves the inductive step and hence Proposition~\ref{p4}.

\medskip

Following the above plan, we choose an eigenvalue $\la_\a$ such that
$\la_\a^{-1}$ is the smallest from among $\la_\b^{-1}, \,\, \b =
1, \dots , p$. Then for every unit vector $X$, the symmetric operator
$\la_\a \L^{-1} - \Phi(X)$ is semidefinite. Indeed, from Lemma~\ref{l5},
$\Phi(X) = N^t N$ for $N \in O_\L = \{N: N \L N^t = \L\}$. Then
$\la_\a \L^{-1} - \Phi(X) = \la_\a \L^{-1} - N^t N = N^t(\la_\a N^{-1t} \L^{-1}
N^{-1} - I_\nu) N = N^t(\la_\a \L^{-1} - I_\nu) N$. The operator
$\la_\a \L^{-1} - I_\nu$ is diagonal in the basis $e_1, \dots, e_\nu$,
with diagonal entries $\tfrac {\la_\a}{\la_\b} - 1 = \la_\a (\la_\b^{-1} -
{\la_\a}^{-1})$. All these numbers have the same sign for $\b \ne \a$.

The $\la_\a$-eigenspace of the Jacobi operator $R_X$ is
$E_{\la_\a}(X) = \{M_Xu : u \in \rnu, \,
R_X (M_Xu) = \la_\a (M_Xu)\}$. In
view of \eqref{e10} the condition $R_X (M_Xu) = \la_\a (M_Xu)$ is equivalent
to $\L \, \Phi(X) u = \la_\a u$, that is, to the fact that $u$ is a
$\la_\a$-eigenvector of the operator $\L \, \Phi(X)$.

Introduce an algebraic set $\mathcal U \subset \Rn \times \rnu$ as follows:
\begin{align*}
\mathcal U &= \{(X, u) \in \Rn \times \rnu : R_X M_Xu = \|X\|^2\la_\a M_Xu \} \\
&=\{(X, u) \in \Rn \times \rnu : \L \, \Phi(X) u = \|X\|^2 \la_\a u \}.
\end{align*}
Let $p_1 : \mathcal U \to \Rn, \,\,p_2 : \mathcal U \to \rnu$ be the projections,
$p_1((X,u)) = X, \linebreak[2] p_2((X,u)) = u$. For every $X \ne 0 \quad p_2 p_1^{-1}(X)$
is the set of vectors $u \in \rnu$ satisfying $\L \Phi(X) u = \|X\|^2 \la_\a u$,
that is, a linear space of dimension $m_\a$, the multiplicity of the
eigenvalue $\la_\a$. So $\dim \mathcal U = n + m_\a$.
It appears that, with our choice of $\la_\a$,
for every $u \in \rnu$ the subset $p_1 p_2^{-1}(u) \subset \Rn$ is also a
linear subspace.

\begin{lemma} \label{l6}
Let $\la_\a$ be the eigenvalue such that
$\la_\a^{-1} = {\mathrm {min}} \{\la_1^{-1}, \dots , \la_p^{-1}\}$. Then for
every $u \in \rnu$ the set
$$
p_1 p_2^{-1}(u) = \{X \in \Rn : \L \Phi(X) u = \|X\|^2 \la_\a u \}
$$
is a linear subspace.
\end{lemma}

\begin{proof} The claim is trivial if $u = 0$. Let $u \ne 0$. The
set $p_1 p_2^{-1}(u)$ is a cone and so it is sufficient to prove that for any
two unit nonparallel vectors $X, \, Y$, the set $p_1 p_2^{-1}(u)$ contains the
unit circle in the
two-plane $\Span(X, Y)$. Let $Z$ be a unit vector in $\Span (X, Y)$ orthogonal
to $X$, and $\T$ be the angle between $X$ and $Y$, so that
$Y = X \cos \T + Z \sin \T$.

Introduce a unit vector function $X_t = \cos t X + \sin t Z$. Then $X = X_0, \,
Y = X_\T$ and we have $\L \Phi(X) u = \L \Phi(X_\T) u = \la_\a u$. As
$\Phi(X_\T) = \cos^2 \T  \, \Phi(X) + \sin^2 \T \, \Phi(Z)
+ \sin \T \cos \T \, (M_X^tM_Z + M_Z^tM_X)$ we find
\begin{equation}\label{e13}
\L \Phi(Z)u = \la_\a u - \cot \T \, \L \, (M_X^tM_Z + M_Z^tM_X) u.
\end{equation}
Let $u_t$ be a $\la_\a$-eigenvector of the operator $\L \Phi(X_t)$ twice
differentiable at $t = 0$ and such that $u_0 = u$. Denote
$\dot u = du_t/dt_{|t=0}, \,\, \ddot u = d^2u_t/dt^2_{|t=0}$.
Differentiating
\begin{equation}\label{e14}
\begin{split}
\L \Phi(X_t) u_t &=(\cos^2 t \, \L \, \Phi(X) + \sin t \cos t \, \L \, (M_X^tM_Z +
M_Z^tM_X)
\\& + \sin^2 t \, \L \, \Phi(Z))u_t = \la_\a u_t
\end{split}
\end{equation}
at $t= 0$ we get $\L(M_X^tM_Z + M_Z^tM_X) u + \L \Phi(X) \dot u = \la_\a \dot u$
and so
\begin{equation}\label{e15}
(M_X^tM_Z + M_Z^tM_X) u = (\la_\a \L^{-1} - \Phi(X)) \dot u.
\end{equation}
The second derivative of \eqref{e14} at $t = 0$ has the form
$$
(-2 \L \Phi(X) + 2 \L \Phi(Z))u + 2 \L \, (M_X^tM_Z + M_Z^tM_X) \dot u + \L \Phi(X)
\ddot u = \la_\a \ddot u.
$$
Substituting the expression for $\L \Phi(Z)u$ from \eqref{e13} we obtain
$$
2 \L(M_X^tM_Z + M_Z^tM_X) (\dot u - \cot \T \, u)  + \L \Phi(X) \ddot u = \la_\a \ddot u.
$$
Acting on both sides by $\L^{-1}$ and taking the inner product with $u$ we get
$$
2 \<(M_X^tM_Z + M_Z^tM_X)u , (\dot u - \cot \T \, u)\> + \< (\Phi(X) -
\la_\a \L^{-1}) u, \ddot u \> = 0.
$$
Substituting $(M_X^tM_Z + M_Z^tM_X) u$ from \eqref{e15} we obtain
$$
2 \<(\la_\a \L^{-1} - \Phi(X)) \dot u, \dot u \> +
\< (\Phi(X) - \la_\a \L^{-1}) u, \ddot u +2 \cot \T \, \dot u \> = 0.
$$
The second term on the left hand side vanishes since
$\L \Phi(X) u = \la_\a u$, hence we get
$$
\<(\la_\a \L^{-1} - \Phi(X)) \dot u, \dot u \> = 0.
$$
With our choice of $\la_\a$, the symmetric operator $\la_\a \L^{-1} - \Phi(X)$
is semidefinite. It follows that $(\la_\a \L^{-1} - \Phi(X)) \dot u = 0$, which
implies $(M_X^tM_Z + M_Z^tM_X) u = 0$ by \eqref{e15}. Then by \eqref{e13}
$\L \Phi(Z)u = \la_\a u$, and so
\begin{align*}
\L \Phi(X_t) u &= (\cos^2 t \, \L \, \Phi(X) + \sin t \cos t \, \L \, (M_X^tM_Z +
M_Z^tM_X)
\\& + \sin^2 t \, \L \, \Phi(Z))u = \la_\a u
\end{align*}
for all $t \in \mathbb R$. It follows that $X_t = \cos t \, X + \sin t \, Z \in
p_1 p_2^{-1}(u)$ for all $t$, that is, $p_1 p_2^{-1}(u)$ contains the
unit circle in the
two-plane $\Span(X, Y)$.
%\qedsymbol
\end{proof}

\begin{lemma} \label{l7}
The subspace $p_2p_1^{-1}(X) \subset \rnu$ is the same for
all $X \ne 0$. In other words, the $\la_\a$-eigenspace of the operator
$\L \Phi(X)$ does not  depend of the choice of a unit vector $X \in \Rn$.
\end{lemma}

\begin{proof} The proof is based on the dimension count.
For every point $u \in \rnu$ let $d(u)$ be the dimension of the linear
space $p_1p_2^{-1}(u)$.

For every set of $m_\a +1$ linearly independent vectors $u_1, \dots ,
u_{m_\a+1}$ in $\rnu$, we must have
\begin{equation}\label{e16}
d(u_1) + \dots + d(u_{m_\a+1}) \le m_\a n.
\end{equation}
Indeed, if the inequality \eqref{e16} is violated, then the subspaces
$p_1p_2^{-1}(u_1), \dots,
 \linebreak[2]
p_1p_2^{-1}(u_{m_\a +1})$ have a nonzero
intersection in $\Rn$. It follows that for a unit vector $X$ from this
intersection, the equation $\L \Phi(X) u = \la_\a u$ has at least
$m_\a +1$ linearly independent solutions, while the $\la_\a$-eigenspace
of the operator $\L \Phi(X)$ has dimension $m_\a$ by 1. of
Lemma~\ref{l5}.

Now let $u_1, \dots , u_{m_\a}$ be a set of $m_\a$ linearly independent vectors
in $\rnu$ such that $d(u_1) + \dots + d(u_{m_\a})$ takes the maximal
possible value.

Let $S = \Span (u_1, \dots , u_{m_\a}), \quad \dim S = m_\a$. For every
point $u \notin S$, we have $d(u) + \sum_{i=1}^{m_\a} d(u_i) \le  m_\a n$
by \eqref{e16}, and $d(u) \le d(u_i)$ for all $i =1, \dots , m_\a$ by the
construction of the set $u_1, \dots , u_{m_\a}$. So
$d(u) \le \tfrac {m_\a}{m_\a +1} n$.

Now the subset $p_2^{-1}(\rnu \setminus S) \in \mathcal U$ is projected by $p_2$
to a subset of $\rnu \setminus S$, of dimension not greater than $\nu$,
with the fibers being linear subspaces of dimension not greater than
$\tfrac {m_\a}{m_\a +1} n$. If $p_2^{-1}(\rnu \setminus S)$ has a nonempty
interior in $\mathcal U$, then $\tfrac {m_\a}{m_\a +1} n + \nu \ge \dim \mathcal U =
n + m_\a$, and so
$$
n \le (\nu - m_\a)(m_\a +1).
$$
This contradicts to \eqref{e9}, since for $1 \le m_\a \le \nu$ the maximal
value of the right hand side is $\tfrac 14 (\nu + 1)^2$.

It follows that the closed subset $p_2^{-1}(S) \in \mathcal U$ must have a
nonempty interior in $\mathcal U$. In a neighbourhood of a generic point from
the interior of $p_2^{-1}(S)$, the projection $p_2$ is a fibration, so
the subspace $S$ contains an open subset $S'$ such that for all
$u \in S', \quad \dim p_2^{-1}(u) = \dim \mathcal U - \dim S = n$. Then
$p_1p_2^{-1}(u) = \Rn$. So for every $X \in \Rn$, the
$\la_\a \|X\|^2$-eigenspace of the operator $\L\Phi(X)$ contains an open
subset of $S$, hence coincides with $S$.
%\qedsymbol
\end{proof}

In Lemma~\ref{l7}, we constructed an $m_\a$-dimensional subspace $S \subset \rnu$
such that for all $X \in \Rn$ and all $u \in S$
$$
\L \Phi(X) u = \la_\a \|X\|^2 u.
$$

Following our plan, we take a basis $u_1, \dots, u_{m_\a}$ in $S$ and define
linear operators $J_s: \Rn \to \Rn, \quad s = 1, \dots, m_\a$ by
$$
J_s X = M_X u_s.
$$
For every unit vector $X$, the $J_sX$'s span $E_{\la_\a}(X)$, the
$\la_\a$-eigenspace of $R_X$. Then $\<J_sX, X\> = 0$ and so all the operators
$J_s$ are skew-symmetric.

To prove the induction step (and hence to prove the Proposition) it remains
to show that, with an appropriate choice of the basis $u_1, \dots, u_{m_\a}$,
the operators $J_s$ are also orthogonal and satisfy the Hurwitz relations.

\begin{lemma} \label{l8}
There exist a basis $u_1, \dots, u_{m_\a}$ in $S$ such
that the operators $J_s$ defined by $J_sX = M_X u_s$ satisfy
$$
J_s J_q + J_qJ_s = -2 \K_{qs} I_n.
$$
for all $s, q = 1, \dots, m_\a$.
\end{lemma}

\begin{proof} Since the operators $J_s$ are skew-symmetric, the
condition $J_s J_q + J_qJ_s = -2 \K_{qs} I_n$ is equivalent to the fact that
for all $X \in \Rn$,
\begin{equation}\label{e17}
\<J_sX, J_qX\> = \K_{qs} \|X\|^2.
\end{equation}

To construct the required basis we pick an arbitrary unit vector
$X_0$ and find an element $N_0 \in \O_\L$ such that
$N_0^t \Phi(X_0) N_0 = I_\nu$ according to Lemma~\ref{l5}.
For any $u \in S$ we have
$\la_\a u = \L \Phi(X_0) u = \L N_0^{-1t} N_0^{-1} u =
N_0 \L N_0^{-1} u$, since $N_0 \L N_0^t = \L$. It follows that
$\L (N_0^{-1} u) = \la_\a (N_0^{-1} u)$, that is, $N_0^{-1} u$ lies
in the $\la_\a$-eigenspace of $\L$, the
coordinate subspace $\Span (e_{m'+1}, \dots, e_{m' + m_\a}) \subset \rnu, \quad
m' = m_1 + \dots + m_{\a-1}$.  Define
$$
u_s = N_0 e_{m'+s}, \quad s = 1, \dots , m_\a.
$$

The operators $J_s$ constructed from this basis satisfy \eqref{e17}.
Indeed, for any unit vector $X$, we have
\begin{align*}
\<J_sX, J_qX\> &= \<M_Xu_s, M_Xu_q\>
= \<\Phi(X)u_s, u_q\>
= \<\la_\a \L^{-1} u_s, u_q\>
\\& = \la_\a \< \L^{-1} N_0 e_{m'+s}, N_0 e_{m'+q}\>
 = \la_\a \< N_0^t \L^{-1} N_0 e_{m'+s}, e_{m'+q}\>
\\& = \la_\a \< \L^{-1} e_{m'+s}, e_{m'+q}\>
= \<e_{m'+s}, e_{m'+q}\> = \K_{sq}.
\end{align*}
%\qedsymbol
\end{proof}

\begin{comment}

\end{comment}


\begin{thebibliography}{99}

\bibitem{1}
J.F.Adams,
Vector fields on spheres,
Bull. Amer. Math. Soc.,
\textbf{68}(1962), 39 -- 41.

\bibitem{2}
D.V.Alekseevsky,
Riemannian spaces with unusual holonomy groups,
Functional Anal. Appl.,
\textbf{2}(1962), 1 -- 10.

\bibitem{3}
M.F.Atiah, R.Bott, A.Shapiro,
Clifford modules,
Topology,
\textbf{3, suppl.1} (1964), 3 -- 38.

\bibitem{4}
R.Brown, A.Gray,
Riemannian manifolds with holonomy group $Spin(9)$,
Diff. Geom. in honor of K.Yano, Kinokuniya, Tokyo (1972), 41 -- 59.

\bibitem{5}
Q.-S.Chi,
A curvature characterization of certain locally rank-one symmetric spaces,
J. Differ. Geom.,
\textbf{28}(1988), 187 -- 202.

\bibitem{6}
E.Garc{\'\i}a-R{\'\i}o, D.N. Kupeli, R.V{\'a}zguez-Lorenzo,
Osserman manifolds in Semi-Riemannian Geometry,
Lecture Notes in Mathematics,
\textbf{1777}(2002), Springer-Verlag.

\bibitem{7}
P.Gilkey,
Manifolds whose curvature operator has constant eigenvalues at the basepoint,
J. Geom. Anal.,
\textbf{4}(1994), 155 -- 158.

\bibitem{8}
P.Gilkey,
Algebraic curvature tensors which are $p$-Osserman,
Diff. Geom. Appl.,
\textbf{14}(2001), 297 -- 311.

\bibitem{9}
P.Gilkey, A.Swann, L.Vanhecke,
Isoperimetric geodesic spheres and a conjecture of Osserman concerning the Jacobi
operator,
Quart. J. Math. Oxford (2),
\textbf{46}(1995), 299 -- 320.


\bibitem{10}
D.Husemoller, Fiber bundles
(1975), Springer-Verlag.

\bibitem{11}
Y.Nikolayevsky,
Osserman manifolds and Clifford structures,
Houston J. Math., {\it to appear}.

\bibitem{12}
Y.Nikolayevsky,
Two theorems on Osserman manifolds,
Diff. Geom. Appl., {\it to appear}.

\bibitem{13}
R.Osserman,
Curvature in the eighties,
Amer. Math. Monthly,
\textbf{97}(1990), 731 -- 756.

\bibitem{14}
Z. Raki\`c,
On duality principle in Osserman manifolds,
Linear Alg. Appl.,
\textbf{296} (1999), 183--189.

\bibitem{15}
W.Thurston,
The theory of foliations of codimension greater than one,
Comment. Math. Helv.,
\textbf{49}(1974), 214 -- 231.

\end{thebibliography}
\end{document}